\documentclass[12pt]{article}
\usepackage{amsmath}
\usepackage[psamsfonts]{amssymb}

\newtheorem{theorem}{\bf Theorem}
\newtheorem{lemma}{\bf Lemma}

\newtheorem{proposition}{\bf Proposition}
\newtheorem{corollary}{\bf Corollary}

\def\field{{\mathbb K}}
\def\C{{\mathbb C}}

\def\0{{\mathbb 0}}

\def\Q{{\mathbb Q}}

\def\Z{{\mathbb Z}}
\def\V{{\mathbb V}}

\title{Integrality and rigidity for postcritically finite polynomials}
\author{Adam Lawrence Epstein}

\begin{document}

\maketitle
\section{Overview}
Parameter spaces of algebraic dynamical systems are stratified according to various natural conditions, such as critical orbit relations, which may be imposed or broken. We consider monic polynomials of degree $d\geq 2$, over fields $\field$ of characteristic zero. Such a polynomial $F$ is determined by a list of critical points ${\bf a}=(a_1,\ldots,a_{d-1})$  and the image $b=F(\bar{\bf a})$ of their barycenter
$$\bar{\bf a}=\frac{1}{d-1}\sum\limits_{i=1}^{d-1}a_i.$$ This prescription yields a family $F_{\bar{\bf a},b}$ parametrized by $\field^{d-1}\times \field$. 

The iterates  of $F$ are defined inductively as 
$${F}^0(z) =  z\;\;\mbox{ and }\;\;{F}^{n+1}(z) = F({F}^n(z)).$$ The {\em forward orbit} of $\zeta\in \field$ is the sequence $F^n(\zeta)$, which is finite if and only if $\zeta$ is preperiodic. More generally, if $\field$ is equipped with a valuation $\nu$, we say that $\zeta$ has {\em bounded} forward orbit if the sequence $\nu(F^n(\zeta))$ is bounded from below. We say $F$ is {\em postcritically finite} (respectively, {\em postcritically bounded}) if the forward orbit of every critical point is finite (respectively, bounded).

An element $\zeta$ of a valued field $\field$ is $\nu$-{\em integral} if $\nu(\zeta)\geq 0$; a $\nu$-integral vector is one whose components are all $\nu$-integral. For prime $p$, we denote by $\nu_p$ the $p$-adic valuation on $\Q$: that is, $\nu_p\left(p^{e}\cdot\frac{x}{y}\right)=e$ if neither integer $x,y$ is a multiple of $p$. An algebraic number is $\nu$-integral for one valuation $\nu$ extending $\nu_p$ if and only if it is for any other, namely if it is the root of a monic polynomial whose coefficients are $\nu_p$-integral rational numbers; thus, we may refer to such an algebraic number is $p$-integral.

\pagebreak
\begin{theorem} \label{integral}
Let $\field$ be a field extending $\Q$, equipped with a valuation $\nu$ extending $\nu_p$. Assume that the degree of $F_{{\bf a},b}$ is a power of $p$, and suppose that $\bar{\bf a}$ is $\nu$-integral. Then $F_{{\bf a},b}$ is postcritically bounded if and only if $\bf a$ and $b$ are $\nu$-integral.
\end{theorem}

The result applies in particular to postcritically finite polynomials, and in this setting sharper results are sometimes possible, at least for certain families. Our interest here is the application of Theorem \ref{integral} to questions of rigidity. It is appropriate to work in the moduli space of translation conjugacy classes, or more concretely, the subspace $\bar{\bf a}=0$ corresponding to {\em centered} polynomials. 
Given $(d-1)$-tuples of integers ${\bf m}=(m_1,\ldots,m_{d-1})$ and ${\bf n}=(n_1,\ldots,n_{d-1})$ with $0\leq m_i<n_i$, consider the vanishing loci $$V_i^{\bf m,n}={\V}\left(
\bar{\bf a},
F^{m_i}_{{\bf a},b}(a_i)-{F}^{n_i}_{{\bf a},b}(a_i)\right)$$ and their intersection 
$V^{\bf m,n}=\bigcap\limits_{i=1}^{d-1}V_i^{\bf m,n}$. 

\begin{corollary} \label{pcf}
Assume that $d$ is a prime power. Then ${V}^{\bf m,n}$ is finite.
\end{corollary}
 
\noindent{\bf Proof: }
Over an algebraically closed field, an algebraic set which is bounded relative to a nontrivial absolute value is necessarily zero-dimensional, hence finite, since any curve projects to a Zariski dense subset of some coordinate axis. In view of  Theorem~\ref{integral}, the $\overline{\Q}$-points of ${V}^{\bf m,n}$ are $p$-integral, whence uniformly bounded relative to any $p$-adic absolute value. Consequently, there are only finitely many $\overline{\Q}$-points, and furthermore, by the Nullstellensatz, every point is a $\overline{\Q}$-point.  $\Box$

\bigskip

In the special case of periodic critical points, we have the following refinement: 

\begin{corollary} \label{main}
Assume that $d$ is a prime power. Then ${V}^{\bf 0,n}$ consists entirely of simple points: near any point of their common intersection, the loci $V_i^{\bf 0,n}$ are smooth, reduced hypersurfaces which are pairwise transverse.
\end{corollary}

\noindent The proof of Corollary \ref{main} is rooted in the well-known computation, usually attributed to Gleason, for monic centered quadratic polynomials: see the Appendix. Bobenrieth has carried out an analogous computation in the family \linebreak $z\mapsto 1+\frac{1}{wz^{d}}$ \cite{bob}. The case of cubic polynomials has also been treated by Silverman, but with $3$-integrality obtained differently, through explicit computation of resultants \cite{silnote}.
 
The restriction to prime power degree is not entirely an artifact of our arithmetic inexpertise. For $n\geq 2$, consider $\psi_n=s_n\circ g_n$, where $s_n(z)=z^n$ and $g_n(z)=\frac{(n+1)z-z^{n+1}}{n}$ (see \cite[Section 4]{pb}). Note that $g_n$ has fixed critical points, at the $n$-th roots of unity, so $\psi_n$ has fixed critical values, at $0$ and $1$. Consequently, any composition $\psi_{n_\ell}\circ\cdots\circ \psi_{n_1}$ is postcritically finite. The leading coefficient of such a composition $\Psi$ is of the form $N^{-1}$ for an integer $N\neq \pm 1$. Conjugation by an appropriate homothety yields a monic centered postcritically finite polynomial having a critical point at a given root $N^{-1/(D-1)}$, where $D$ is the degree of $\Psi$, and such a point cannot be $\nu_p$-integral for any $p$ dividing $N$.  For example, $\psi_2$ rescales to $z^6-(2^{1/5}3)z^4 + (2^{-8/5}3^2) z^2$ with a critical point at $2^{-2/5}$: this critical point is not $2$-integral, and yet $2$ divides the degree $6$. In this case there is another prime dividing $6$, namely $3$, and all of the critical points are $3$-integral. However, for the degree $72$ polynomial $\psi_2\circ \psi_3$ (respectively $\psi_3\circ \psi_2$), the critical point $2^{-{2/71}}3^{-18/71}$ (respectively $-2^{-24/71}3^{-3/71}$) is not $p$-integral for either of the primes $p=2,3$ dividing $72$. The appropriate extension of Theorem \ref{integral} remains a mystery.

On the other hand, Corollary \ref{pcf} (without degree restrictions) is well-known through the use of complex analytic methods, as the analogous corollary of the boundedness of the {\em connectedness locus}: the set of parameters corresponding to polynomials which are postcritically bounded relative to the usual archimedean absolute value, or more customarily, whose Julia set over $\C$ is connected. Theorem~\ref{integral} is a nonarchimedean version of this well-known boundedness principle, and our proof is similar in spirit. We imagine that the corresponding archimedean computation is folklore: the now classical argument in \cite{BH} relies on basic estimates from univalent function theory.

Corollary \ref{main} (without degree restrictions) is also understood complex analytically, but the explanation is conceptually deeper. The relevant transversality is deduced in various ways. For example, \cite{milnor} classifies the biholomorphism types of hyperbolic components, and verifies transversality by inspection of the models. The discussion in \cite{E} from first principles of Teichm\"uller theory additionally yields appropriate transversality assertions in connection with preperiodic critical points. Both treatments are conceptually related to a fundamental rigidity principle due to Thurston (see \cite{DH2} for discussion and proof of the Existence and Uniqueness Theorems). We remark that the algebraic content of Thurston Rigidity for polynomials is expressed by (unrestricted) Corollary \ref{pcf}, and that the proof of Thurston Rigidity rests on an infinitesimal rigidity principle whose algebraic content, for polynomials with periodic critical points, is expressed by (unrestricted) Corollary~\ref{main}.

\subsection*{Acknowledgments}
We thank Bjorn Poonen and Joe Silverman for enlightening discussions of related matters. We thank Xander Faber for organizing the May, 2010 workshop on Moduli Spaces and the Arithmetic of Dynamical Systems, at the Bellairs Research Institute in Barbados, where these discussions arose. We thank Xavier Buff for suggesting  \cite[Section 4]{pb} as a source of counterexamples.

\section{Integrality}
Consider the polynomials
$$F_{{\bf a},b}(z)=z^d+\sum_{k=1}^{d-1}(-1)^{d-k}\frac{d}{k}\sigma_{d-k}z^k+b-\bar{\bf a}^d-\sum_{k=1}^{d-1}(-1)^{d-k}\frac{d}{k}\sigma_{d-k}{\bar{\bf a}}^k$$
where
$$\sigma_k=\sum\limits_{1\leq i_1<\cdots<i_k\leq d-1}\prod\limits _{j=1}^k a_{i_j}$$
are the elementary symmetric functions.
Observe that $b$ is the image of the barycenter $\bar{\bf a}=\frac{1}{d-1}\sum\limits_{i=1}^{d-1}a_i$ of the (labeled) critical points $a_1,\ldots,a_{d-1}$. We denote by $c$ the barycenter $$\frac{1}{d-1}\sum_{i=1}^{d-1} F_{{\bf a},b}(a_i)$$ of the critical values.
We set ${\bf a}^*=(a_1^*,\ldots,a_{d-1}^*)$, where $a_i^*=a_i-\bar{\bf a}$, and $b^*=b-\bar{\bf a}$, so that $F_{{\bf a},b}(z+\bar{{\bf a}})=F_{{\bf a}^*,b^*}(z)+\bar{\bf a}$. Note that $c=b+\Phi({\bf a}^*)$
where $$\Phi({\bf a})=\frac{1}{d-1}\sum_{i=1}^{d-1}\left(a_i^d+\sum_{k=1}^{d-1}(-1)^{d-k}\frac{d}{k}\sigma_{d-k}a_i^{k}\right)
-\bar{\bf a}^d-\sum_{k=1}^{d-1}(-1)^{d-k}\frac{d}{k}\sigma_{d-k}{\bar{\bf a}}^k$$
is homogeneous of degree $d$, and translation invariant: $\Phi({\bf a})=\Phi({\bf a}^*)$.

The fact that the coefficients of $F_{{\bf a},b}$ lie in $\Q$, but not necessarily in $\Z$, is a source of complication. However, under favorable conditions, these coefficients do belong to appropriate local rings $\Z_{(p)}$. 
\begin{lemma} \label{pp}
Assume that $d$ is a power of a prime $p$. Then:
\begin{itemize}
\item $F_{{\bf a},b}(z)\equiv z^d+b-\bar{\bf a}^d\;({\rm mod}\;p)$ in $\Z_{(p)}[a_1\ldots, a_{d-1},b][z]$,
\item $\Phi({\bf a})\equiv 0\;({\rm mod}\;p)$ in $\Z_{(p)}[a_1,\ldots a_{d-1}]$.
\end{itemize}
\end{lemma}

\noindent{\bf Proof: } Note that if $p|d$ then $\frac{1}{d-1}\in\Z_{(p)}$ so $\bar{\bf a}\in\Z[a_1,\ldots,a_{d-1}]$. Moreover, if
$d$ is a power of $p$ then $\nu_p\left(\frac{d}{k}\right)\geq 1$ for $1\leq k\leq d-1$, and furthermore $$\Phi({\bf a})\;\equiv\;\frac{1}{d-1}\sum_{i=1}^{d-1} a_i^d  -\left(\frac{1}{d-1}\sum_{i=1}^{d-1}a_i\right)^{d}\;({\rm mod}\;p).$$
Now 
$$\left(\sum_{i=1}^{d-1}a_i\right)^d= \sum_{(m_1 \ldots m_{d-1})}{d\choose m_1 \ldots m_{d-1}}\prod_{i=1}^{d-1} a_i^{m_i}$$ summed over $(d-1)$-tuples of integers $m_i\geq 0$ with $\sum\limits_{i=1}^{d-1}m_i=d$. Since
$$\nu_p(n!)=\sum_{e=1}^\infty\left\lfloor\frac{n}{p^e}\right\rfloor$$
for any $n$, we have
$$\nu_p\left(n\choose m_1 \ldots m_{d-1}\right)=\nu_p(n!)-\sum_{i=1}^{d-1}\nu_p(m_i!)=\sum_{e= 1}^\infty\left(\left\lfloor\frac{n}{p^e}\right\rfloor-\sum_{i=1}^{d-1}\left\lfloor\frac{m_i}{p^e}\right\rfloor\right)$$
where $\left\lfloor\frac{n}{p^e}\right\rfloor\geq\sum\limits_{i=1}^{d-1} \left\lfloor\frac{m_i}{p^e}\right\rfloor$ for every $e\geq 1$. Thus, if $d=p^e$ and all $m_i<d$ then $$\nu_p\left(d\choose m_1 \ldots m_{d-1}\right)\geq \left\lfloor\frac{p^e}{p^e}\right\rfloor-\sum_{i=1}^{d-1}\left\lfloor\frac{m_i}{p^e}\right\rfloor=1,$$ 
hence
$\left(\sum\limits_{i=1}^{d-1}a_i\right)^d\equiv\sum\limits_{i=1}^{d-1} a_i^d \;({\rm mod}\;p)$. Since $d|((d-1)^{d-1}-1)$, it follows that $$\Phi({\bf a})\equiv((d-1)^{d-1}-1)\bar{\bf a}^d\equiv 0\;({\rm mod}\;p).$$
$\Box$
\pagebreak

\bigskip
Let $\nu$ be a valuation extending $\nu_p$, where $p|d$.
For $({\bf a},b)\in \field^{d-1}\times \field$, set
$$\mu_{{\bf a},b}=\inf\left\{\nu(\zeta-\bar{\bf a}): \mbox{ the forward orbit } F_{{\bf a},b}^n(\zeta) \mbox{ is bounded}\right\}.$$
Since $\nu(F^n_{{\bf a},b}(\zeta)-F^n_{{\bf a}^*,b^*}(\zeta^*))=\nu(\bar{\bf a})$ is constant, for any $\zeta\in \field$ and $\zeta^*=\zeta-\bar{\bf a}$, the forward orbit $F^n_{{\bf a},b}(\zeta)$ is bounded if and only if the forward orbit $F^n_{{\bf a}^*,b^*}(\zeta^*)$ is bounded. It follows that $\mu_{{\bf a,}b}=\mu_{{\bf a}^*,b^*}$ and that $F_{{\bf a},b}$ is postcritically bounded if and only if $F_{{\bf a}^*,b^*}$ is postcritically bounded.

\begin{lemma} \label{bddorbit}
For any $({\bf a},b)\in \field^{d-1}\times \field$, we have
$$\mu_{{\bf a},b}\geq \min\left(\alpha+\epsilon,\frac{\beta}{d},0\right)$$ where $\alpha=\min\limits_{1\leq i\leq d-1}\nu(a_i-\bar{\bf a})$
and $\beta=\nu(b-\bar{\bf a})$, and $\epsilon=\min\limits_{1\leq k\leq d-1}\frac{\nu\left(\frac{d}{k}\right)}{d-k}$.
\end{lemma}

\noindent {\bf Proof: } Since the statement is translation invariant, we may assume without loss of generality that $\bar{\bf a}=0$, so $\alpha=\min\limits_{1\leq i\leq d-1}\nu(a_i)$ and therefore
$$\nu\left(\frac{d}{k}\right)+(d-k)\alpha+\nu(\zeta^{k})\leq\nu\left(\frac{d}{k}\sigma_{d-k} \zeta^{k}\right).$$ 
Thus, if $\nu(\zeta)<\alpha+\frac{\nu\left(\frac{d}{k}\right)}{d-k}$ then $\nu(\zeta^d)=\nu(\zeta^{d-k})+\nu(\zeta^{k})<\nu\left(\frac{d}{k}\sigma_{d-k} \zeta^{k}\right)$,
so if $\nu(\zeta)<\alpha+\epsilon$ then $\nu(\zeta^d)<\nu\left(\frac{d}{k}\sigma_{d-k} \zeta^{k}\right)$
for $1\leq k\leq d-1$, whence if $\nu(\zeta)<\min \left(\alpha+\epsilon,\frac{\beta}{d}\right)$ then $\nu(F_{{\bf a},b}(\zeta))=\nu(\zeta^d)=d\nu(\zeta)$. Consequently, if $\nu(\zeta)<\min\left(\alpha+\epsilon,\frac{\beta}{d},0\right)$ then $\nu(F_{{\bf a},b}(\zeta))=d\nu(\zeta)<\nu(\zeta)<\min\left(\alpha+\epsilon,\frac{\beta}{d},0\right)$, whence $\nu(F^n_{{\bf a},b}(\zeta))=d^n\nu(\zeta)\to -\infty$. $\Box$

\bigskip
\noindent {\bf Proof of Theorem \ref{integral}: } As above, we may assume without loss of generality that $\bar{\bf a}=0$.
By Lemma \ref{pp}, since $d$ is a power of $p$ we have $\epsilon=\frac{1}{d-1}>0$
and moreover $\nu(c-b)=\nu(\Phi({\bf a}))\geq 1+d\alpha$, since $\Phi$ is homogeneous of degree $d$. 
Now if $F_{{\bf a},b}$ is postcritically bounded then
$\alpha\geq \mu_{{\bf a},b}$, hence $\mu_{{\bf a},b}\geq\min\left(\frac{\beta}{d},0\right)$, and
$\nu(c)\geq \min\limits_{1\leq k\leq d-1}\nu(F_{{\bf a},b}(a_i))\geq\mu_{{\bf a},b}$. Thus,
$$\beta\geq\min(\nu(c),\nu(c-b))\geq\min(\mu_{{\bf a},b},1+d\mu_{{\bf a},b})\geq\min\left(\frac{\beta}{d},0,1+\beta,1 \right)$$
so $\beta\geq\min\left(\frac{\beta}{d},0\right)$, hence $\beta\geq 0$, whence $\alpha\geq\mu_{{\bf a},b}\geq\frac{\beta}{d}\geq 0$: that is, ${\bf a}$ and $b$ are $\nu$-integral.
Conversely, if $\bf a$ and $b$ are $\nu$-integral then since $\epsilon\geq 0$ the coefficients of $F_{{\bf a},b}$ are $\nu$-integral, so the postcritical points are $\nu$-integral, whence $F_{{\bf a},b}$ is postcritically bounded.  $\Box$

\bigskip
The well-known considerations of Lemma \ref{bddorbit} establish the existence and basic properties of the {\em local canonical height} functions
$$h_{{\bf a},b}(\zeta)=\lim_{n\to\infty}\frac{1}{d^n}\max(-\nu(F^n_{{\bf a},b}(\zeta),0))$$
as discussed in the arithmetic dynamics literature \cite{sil}; in these terms, $F_{{\bf a},b}$ is postcritically bounded if and only if $H_{{\bf a},b}=0$, where $$H_{{\bf a},b}=\max_{1\leq i\leq d-1}h_{{\bf a},b}(a_i).$$ These quantities are evidently the nonarchimedean analogues of the Green's functions used in the classical proof of the boundedness of the connected locus in the archimedean case.

\section{Simplicity}

Recall that the loci ${V}^{\bf m,n}$ were defined in terms of the parameter subspace of centered polynomials $F_{{\bf a},b}$. Here it will be convenient to work with the variant family
$${\cal F}_{\bf a}(z)=z^d+\sum_{k=1}^{d-1}(-1)^{d-k}\frac{d}{k}\sigma_{d-k}z^k$$
normalized to fix $0$. Since $d$ is a power of $p$, it follows from Lemma \ref{pp} 
that ${\cal F}_{{\bf a}}(z)\equiv z^d\;({\rm mod}\;p)$ in $\Z_{(p)}\left[a_1,\ldots,a_{d-1},b\right][z]$.  Consequently,
$$\frac{\partial {\cal F}_{{\bf a}}(w)}{\partial a_i}=dw^{d-1}\frac{\partial w}{\partial a_i}\equiv 0 \;({\rm mod}\;p)$$
for $1\leq i\leq d-1$ and any $w\in \Z_{(p)}\left[a_1,\ldots,a_{d-1},b\right]$. 

Note that translation by $\bar{\bf a}$ conjugates ${\cal F}_{\bf a}$ to $F_{{\bf a}^*,b^*}$ where 
$$({\bf a}^*,b^*)=\left(a_1-\bar{\bf a},\ldots,a_{d-1}-\bar{\bf a},{\cal F}_{{\bf a}}(\bar{\bf a})-\bar{\bf a}\right).$$
Conversely, if $\bar{\bf a}^\star=0$ then any fixed point $\xi$ of $F_{{\bf a}^\star,b^\star}$ we have $({\bf a}^\xi)^*={\bf a}^\star$
for ${\bf a}^\xi=(a^\star_1-\xi,\ldots,a^\star_{d-1}-\xi)$. Observe that the map ${\bf a}\mapsto ({\bf a}^*,b^*)=\Lambda({\bf a})$ sends each locus ${\cal V}_i^{\bf m,n}={\V}({\cal F}_{\bf a}(a_i)-a_i)$ onto the corresponding locus $V_i^{\bf m,n}$. Moreover, $\Lambda$ respects $p$-integrality:
if ${\bf a}$ is $p$-integral then $\bar{\bf a}$ is $p$-integral, hence $({\bf a}^*,b^*)$ is also $p$-integral, while if $({\bf a}^\star,b^\star)$ is $p$-integral then, by monicity, any fixed point $\xi$ of $F_{{\bf a}^\star,b^\star}$ is $p$-integral, whence ${\bf a}^\xi$ is $p$-integral. 

\bigskip
\noindent {\bf Proof of Corollary \ref{main}: } We claim that $\Lambda$
is nonsingular at every $p$-integral ${\bf a}$. Indeed, the derivative of the composition ${\bf a}\mapsto({\bf a}^*,b^*)\mapsto(a_1,\ldots,a_{d-2},b)$ is given by the $(d-1)\times(d-1)$ matrix
$$\left(\begin{array}{rrrrr}
1-\frac{1}{d-1} & -\frac{1}{d-1} & \cdots & -\frac{1}{d-1} & -\frac{1}{d-1}\\
-\frac{1}{d-1}  & 1-\frac{1}{d-1} &\cdots & -\frac{1}{d-1}& -\frac{1}{d-1}\\
\cdots &\cdots &\cdots & \cdots &\cdots \\
 -\frac{1}{d-1}& -\frac{1}{d-1}&\cdots & 1-\frac{1}{d-1}& -\frac{1}{d-1}  \\
\frac{\partial {\cal F}_{{\bf a}}(\bar{\bf a})}{\partial a_1} -\frac{1}{d-1}& \frac{\partial {\cal F}_{{\bf a}}(\bar{\bf a})}{\partial a_2} -\frac{1}{d-1}&\cdots & \frac{\partial {\cal F}_{{\bf a}}(\bar{\bf a})}{\partial a_{d-2}}-\frac{1}{d-1} &\frac{\partial {\cal F}_{{\bf a}}(\bar{\bf a})}{\partial a_{d-1}} -\frac{1}{d-1}
\end{array}\right),$$
and since $\bar{\bf a}\in\Z_{(p)}[a_1,\ldots,a_{d-1},b]$, this matrix is congruent $({\rm mod}\;p)$ to
$$\left(\begin{array}{rrrrr}
0 & -1 & \cdots & -1 & -1\\
-1  &0 &\cdots & -1& -1\\
\cdots &\cdots &\cdots & \cdots &\cdots \\
 -1& -1&\cdots & 0& -1 \\
-1& -1&\cdots & -1&-1
\end{array}\right)$$
which has determinant $-1\not \equiv 0$. In view of Theorem \ref{integral}, every point of ${V}^{\bf 0,n}$ is $p$-integral. It follows from the remarks above that $\Lambda$ yields local isomorphisms, respecting integrality, between neighborhoods of points in ${\cal V}^{\bf 0,n}$ and neighborhoods of the corresponding points in ${V}^{\bf 0,n}$. Moreover, proving simplicity for ${\cal V}^{\bf 0,n}$ amounts to observing the invertibility of ${\bf I}-{\bf F}$, where $\bf I$ is the $(d-1)\times(d-1)$ identity matrix, and where
$${\bf F}=\left(
\begin{array}{ccc}
\frac{\partial {\cal F}^{n_1}_{{\bf a}}(a_1)}{\partial a_1} &\cdots & \frac{\partial {\cal F}^{n_1}_{{\bf a}}(a_1)}{\partial a_{d-1}}\\
\cdots & \cdots & \cdots \\
\frac{\partial {\cal F}^{n_{d-1}}_{{\bf a}}( a_{d-1})}{\partial a_1} &\cdots & \frac{\partial {\cal F}^{n_{d-1}}_{{\bf a}}(a_{d-1})}{\partial a_{d-1}}
\end{array}
\right)\equiv 0\;({\rm mod}\;p)$$ since ${\cal F}^{n_{i}-1}_{{\bf a}}( a_{i})\in \Z_{(p)}\left[a_1,\ldots,a_{d-1}\right]$ for $1\leq i\leq d-1$. $\Box$

\appendix
\section*{Appendix (joint with Bjorn Poonen)}

Consider the family of monic centered unicritical polynomials $F_{{\bf 0},b}(z)=z^d+b$, where $d$ is any integer greater than 1. For $n\geq 0$, we set $\Gamma_n=F_{{\bf 0},b}^n(0)\in\Z[b]$. Note that $\Gamma_0=0$, and that $\Gamma_n$ is monic of degree $2^{n-1}$, for $n\geq 1$. Thus, for $N>n\geq 0$ the polynomials $\Gamma_N-\Gamma_n\in\Z[b]$ are monic, whence their zeros in $\overline{\Q}$ are algebraic integers, in accordance with Theorem~\ref{integral}.

Here we extend Corollary \ref{main} to all postcritically finite parameters in these families. The case of periodic critical point is straightforward. For $p=d=2$, the following observation, already contained in Corollary \ref{main}, is due to Andrew Gleason (see \cite[Lemma 19.1]{DH}), and independently to Allen Adler: 

\begin{proposition}\label{g}
For $N\geq 1$, all zeros of $\Gamma_N$ are simple.
\end{proposition}

\noindent {\bf Proof: } By definition, $\Gamma_{N} = (\Gamma_{N-1})^d+b$, so $\Gamma_{N}' = d (\Gamma_{N-1})^{d-1} \Gamma_{N-1}' + 1$, and thus $\Gamma_{N}'\equiv 1 \;({\rm mod}\;p)$ for any prime $p|d$. It follows that $\Gamma_N'(b)\neq 0$ for every algebraic integer $b$, in particular, for every zero of $\Gamma_N$. $\Box$

\bigskip
The case of strictly preperiodic critical point is more subtle. For $d=2$ this is treated in \cite[Lemma 1, page 333]{DH3}, but the proof is incomplete when $n=2$, since the discussion presumes that $\Gamma_{n-2}' \equiv 1 \;({\rm mod}\;2)$.

\begin{theorem}
For $N>n\geq 1$, every multiple zero of $\Gamma_{N}-\Gamma_n$ is a zero of $\Gamma_{N-1}-\Gamma_{n-1}$. 
\end{theorem}

\noindent {\bf Proof: } Observe that $\Gamma_{N}-\Gamma_n=(\Gamma_{N-1})^d-(\Gamma_{n-1})^d=\prod_\omega\Delta_{N,n}^\omega$ where $\Delta_{N,n}^\omega=\Gamma_{N-1}-\omega \Gamma_{n-1}$ and where the product is over the $d$-th roots of unity. If $\omega_1\neq\omega_2$ then any common zero of $\Delta_{N,n}^{\omega_1}$ and $\Delta_{N,n}^{\omega_2}$ is also a zero of $\Gamma_{N-1}$ and $\Gamma_{n-1}$, whence a zero of $\Delta^\omega_{N,n}$ for every $\omega$, in particular for $\omega=1$. Thus, it suffices to show that if $\omega\neq 1$ then all zeros of $\Delta^\omega_{N,n}$ are simple. Since every zero of $\Gamma_{N}-\Gamma_1=(\Gamma_{N-1}-\Gamma_0)^d$ is a zero of $\Gamma_{N-1}-\Gamma_0$, we may assume without loss of generality that $n\geq 2$, whence 
$\Gamma_{N-1}'\equiv 1 \equiv \Gamma_{n-1}\;({\rm mod}\;p)$ for any prime $p|d$. It follows that if $\nu$ is a valuation extending $\nu_p$ then $\nu({\Delta^\omega}'_{N,n}(b)-(1-\omega))\geq 1$ for every algebraic integer $b$; this holds in particular for the zeros of $\Delta^\omega_{N,n}$, so if $\nu(1-\omega)<1$ then these zeros are simple.

Each $\omega$ is a primitive $m$-th root of unity for some $1\neq m|d$. If $m=p^e$ for some $e\geq 1$ then $\nu(1-\omega)=(p-1)p^{e-1}$ for any $\nu$ extending $\nu_p$, while if $m$ is not a prime power then $1-\omega$ is a unit so $\nu(1-\omega)=0$ for any valuation $\nu$: for details, see \cite[page 73]{lang}. Since $(p-1)p^{e-1}>1$ except when $p=2$ and $e=1$, it follows that if $m\neq 2$ then there exists $p|d$ such that $\nu(1-\omega)<1$ for any $\nu$ extending $\nu_p$. Furthermore, if $m=2$ then $\omega=-1$, so if $d$ has an odd prime factor $p$ then $\nu(1-\omega)=\nu_p(2)=0$. These considerations establish simplicity in all cases except when $\omega=-1$ and $d$ is a power of 2.

Suppose finally that $\omega=-1$ and $d=2^e$ for some $e\geq 1$. In this case, ${\Delta^\omega}'_{N,n}=\Gamma_{N-1}'+\Gamma_{n-1}'=d\Lambda+2$  where $\Lambda=(\Gamma_{N-2})^{d-1}\Gamma_{N-2}'+(\Gamma_{n-2})^{d-1}\Gamma_{n-2}'$. Consequently, it suffices to show that $\Delta_{N,n}^\omega(b)=0$ implies $\nu(\Lambda(b))>0$, since then $\nu(d\Lambda(b))> e \geq 1=\nu(2)$. Observe that if $\nu(x),\nu(y)\geq 0$ \linebreak then $\nu((x+y)-(x-y))=\nu(2y)\geq 1$, so if $\nu(x-y)=0$ then $\nu(x+y)=0$, hence $\nu(x^2-y^2)=0$, and thus $\nu(x^{2^k}-y^{2^k})=0$ for $k\geq 0$; in particular, since $F_{{\bf 0},b}(x)-F_{{\bf 0},b}(y)=x^{2^e}-y^{2^e}$, it follows that $\nu(F_{{\bf 0},b}(x)-F_{{\bf 0},b}(y))>0$ implies $\nu(x-y)>0$. Now if $\Delta^\omega_{N,n}(b)=0$ then $F_{{\bf 0},b}^2(\Gamma_{N-2}(b))=F_{{\bf 0},b}^2(\Gamma_{n-2}(b))$, so $\nu(\Gamma_{N-1}(b)-\Gamma_{n-1}(b))>0$ and thus $\nu(\Gamma_{N-2}(b)-\Gamma_{n-2}(b))>0$. If $n>2$ then $\nu(\Gamma_{N-2}'(b)-1),\nu(\Gamma_{n-2}'(b)-1)\geq 1$, and thus $\nu(\Gamma_{N-2}'(b)+\Gamma_{n-2}'(b))\geq 1$; since
$\nu(\Gamma_{N-2}(b)^{d-1}-\Gamma_{n-2}(b)^{d-1})>0$ and $\nu(\Gamma_{n-2}(b))\geq 0$, we have 
$\nu(\Lambda(b))>0$. If $n=2$ then $\Gamma_{n-2}(b)=0$ so $\nu(\Gamma_{N-2}(b))>0$; since $\nu(\Gamma_{N-2}'(b))\geq 0$, it follows that $\nu(\Lambda(b))>0$ in this case as well.  $\Box$

\end{document}